\documentclass[12pt]{article}

\usepackage{graphics,amsmath,makeidx}
\usepackage{amssymb}  
\usepackage{amsfonts}
\usepackage{graphicx}

\newtheorem{que}{Question}[section]
\newtheorem{thm}{Theorem}[section]
\newtheorem{lem}{Lemma}[section]
\newtheorem{dfn}{Definition}[section]
\newtheorem{ex}{Example}[section]

\newenvironment{pf}{\SK\noindent{\bf Proof.}\small}{\hfill
$\Box$\smallbreak\SK}

\newcommand{\rf}[1]{(\ref{#1})}
\newcommand{\cl}[1]{{\cal #1}}

\newcommand{\tx}[1]{\mbox{\rm{#1}}}

\def\Box{\hbox{\hskip 1pt \vrule width 8pt height 6pt depth 1.5pt
  \hskip 1pt}}

\renewcommand{\a}{\alpha}
\renewcommand{\b}{\beta}
\newcommand{\e}{\epsilon}
\newcommand{\g}{\gamma}

\renewcommand{\d}{\delta}
\newcommand{\D}{\Delta}

\newcommand{\W}{\Omega}

\renewcommand{\t}{\theta}

\newcommand{\s}{\sigma}

\newcommand{\bydef}{\stackrel{\mathrm{def}}{=}}

\newcommand{\EQ}[2]{\begin{equation}\label{#1}#2\end{equation}}

\newcommand{\AR}[2]{\left[\begin{array}{#1}#2\end{array}\right]}

\newcommand{\OR}[2]{\left\{\begin{array}{#1}#2\end{array}\right.}
\newcommand{\NI}{\noindent}
\newcommand{\SK}{\vskip2mm}

\newcommand{\tr}{\mbox{{\rm tr}}}

\newcommand{\RR}{\mathbb{R}}
\newcommand{\CC}{\mathbb{C}}
\newcommand{\TT}{\mathbb{T}}
\newcommand{\DD}{\mathbb{D}}
\newcommand{\ZZ}{\mathbb{Z}}

\newcommand{\THM}[2]{\begin{thm}\label{#1}{\sl #2}\end{thm}}

\title{KYP Lemma for Non-Strict Inequalities\\ and 
the associated Minimax Theorem}

\author{A. Megretski, MIT (ameg@mit.edu) 
}

\begin{document}

\maketitle

\begin{abstract}
\NI Several variations of the classical
Kalman-Yakubovich-Popov Lemma, as well the 
associated minimax theorem are presented.
\end{abstract}

\subsubsection*{Notation and Terminology}
$\ZZ_+$ is the set $\ZZ_+=\{0,1,2,\dots\}$ of all non-negative integers.
$j\RR=\{s\in\CC:\ \tx{Re}(s)=0\}$,
$\CC_+= \{s\in\CC:\ \tx{Re}(s)>0\}$, $\TT=\{z\in\CC:\ |z|=1\}$, and 
$\DD_+=\{z\in\CC:\ |z|>1\}$
are the frequently referenced subsets of the complex plane $\CC$:
the imaginary axis, the open right half plane,
the unit circle, and the outside of the unit circle.
$\CC^{n,m}\supset\RR^{n,m}$ are the sets of $n$-by-$m$
matrices (complex and real), with the usual shortcuts
$\CC^n=\CC^{n,1}$, $\RR^n=\RR^{n,1}$. For $M\in\CC^{k,n}$,
$M'\in\CC^{n,k}$ is the Hermitian conjugate (the result of applying
both transposition and
complex conjugation to $M$), while $\bar M\in\CC^{k,n}$
is the complex conjugate of $M$. For a real vector space $V$,
$V^\sharp$ is the real vector space of all linear functionals
$f:\ V\mapsto\RR$.

\section{The Classical KYP Lemma}
A number of alternative versions of
the KYP Lemma,
a classical result of the linear system theory, 
has been published over the last half century.
The earlier formulations, such as \cite{Yak1}, 
motivated by optimal linear
feedback design applications, related positive definiteness (or
semi-definiteness) of rational matrix-valued functions of a single
complex variable on the real axis or on the unit circle
(the so-called "frequency conditions") to the existence of
"stabilizing" (or "marginally stabilizing")
solutions of the associated Lur'e (algebraic Riccati) equations.
Connections to dynamic programming and first order conditions
of optimality allowed extensions to time-varying and distributed
systems (see, for example, \cite{Yak2,Yak3}).
Some of the more recent versions, such as \cite{Ran1}, 
employ weaker
assumptions to relate the frequency domain inequalities
to feasibility of the semidefinite programs obtained by
replacing the Lur'e or Riccati equations by the corresponding
inequalities. 

It appears that
some useful versions of the KYP Lemma remain unpublished
(or, at least, highly inaccessible). This paper aims to
correct this by presenting 
several (assumedly) missing formulations. 

\subsection{KYP Lemma in Discrete Time}
The classical KYP Lemma setup is
defined by matrices
$A\in\CC^{n,n}$, $B\in\CC^{n,m}$, $Q=Q'\in\CC^{n+m,n+m}$:
$A$ and $B$ are the coefficients of linear transformation 
\[ (x,u)\in\CC^n\times\CC^m\mapsto x_+=Ax+Bu\in\CC^n,\]
and $Q$ is associated with the
Hermitian form $\s:\ \CC^n\times\CC^m\mapsto\RR$:
\EQ{k1}{  \s(x,u)=
\AR{c}{x\\ u}'Q\AR{c}{x\\ u}\ \ \ (x\in\CC^n,\ u\in\CC^m).}

\subsubsection{Stabilizing Completion of Squares in Discrete Time}

This is one of the versions of the KYP Lemma, motivated by
the linear quadratic optimal control design theory.

\THM{thm:kypstabdt}{
Assume that the pair $(A,B)$, where 
$A\in\CC^{n,n}$ and $B\in\CC^{n,m}$, 
is stabilizable, in the sense that
there exists a matrix $K\in\CC^{m,n}$ such that $zI_n-A-BK$ is
invertible for all $z\in\CC$, $|z|\ge1$.
Then for every matrix
$Q=Q'\in\CC^{n+m,n+m}$
(and $\s$ defined in \rf{k1})
the following conditions are equivalent:
\begin{itemize}
\item[(a)] there exist matrices
$P=P'\in\CC^{n,n}$, $C\in\CC^{m,n}$, $D\in\CC^{m,m}$ 
such that  
\EQ{k2}{\s(x,u)+x'Px-(Ax+Bu)'P(Ax+Bu)=|Cx+Du|^2\ \ \ 
\forall\ x\in\CC^n,\ u\in\CC^m,}
\EQ{k3}{ 
\det\AR{cc}{\lambda A-I_n & \lambda B\\ C & D}\neq0\ \ 
\forall\ |\lambda|<1;}
\item[(b)]  the matrix 
\EQ{k4}{\Pi(z)=
\AR{c}{(zI_n-A)^{-1}B\\ I_m}'Q\AR{c}{(zI_n-A)^{-1}B\\ I_m},}
defined for $z\not\in\Lambda(A)=\{z\in\CC:\ \det(zI_n-A)=0\}$,
is positive definite for all $z\in\TT$
except, possibly,
a finite subset.
\end{itemize}
Moreover, when matrices $A,B,Q$ in (b) are real,
the corresponding matrices $P,C,D$ from (a) can be chosen to be
real as well.
}

A proof of Theorem \ref{thm:kypstabdt} 
(as well as any other proof provided in this paper)
can be found in the Appendix
section.

We will refer to Theorem~\ref{thm:kypstabdt} as the
"stabilizing completion of squares"
version of the KYP Lemma, because the right side of \rf{k2} can
be viewed of a "complete square" quadratic form, and 
\rf{k3} guarantees that the matrix $A-BD^{-1}C$ is well defined and
"marginally stable" (has no eigenvalues $z$ with $|z|>1$).

\subsubsection{Application: Optimal Program Control}

The "stabilizing completion of squares" was originally motivated by
an "abstract" optimal control question of finding the
maximal lower bound of the functional
\EQ{k7}{ \Phi(x(\cdot),u(\cdot))=\sum_{t=0}^\infty \s(x(t),u(t))\to\inf}
subject to linear equations 
\EQ{k7A}{x(t+1)=Ax(t)+Bu(t),\ x(0)=a,}
and the "finite energy" constraint
\EQ{k7B}{\sum_{t=0}^\infty\{|x(t)|^2+|u(t)|^2\}<\infty,}
where $A,B,Q,a$ are fixed, and
$x:\ \ZZ_+\mapsto\CC^n$, $u:\ \ZZ_+\mapsto\CC^m$
are infinite dimensional decision variables. The following statement,
which follows directly from Theorem~\ref{thm:kypstabdt},
explains the relation between the optimization setup 
\rf{k7}-\rf{k7B} and
Theorem~\ref{thm:kypstabdt}.

\THM{thm:lqrdt}{
If the stabilizability
assumption as well as
conditions (a),(b) from Theorem~\ref{thm:kypstabdt} are satisfied then
the infimum in \rf{k7}-\rf{k7B} 
equals $-a'Pa$, and the sum
\[\sum_{t=0}^\infty|Cx_i(t)+Du_i(t)|^2\] 
converges
to zero if and only if $\Phi(x_i(\cdot),u_i(\cdot))$ converges to $-a'Pa$.
}

\subsubsection{Strict Linear Matrix Inequalities}

In many applications, the "stabilization" constraint is irrelevant,
which motivates the following "strict linear matrix inequality (LMI)"
version of the KYP Lemma. 

\THM{thm:kyplmisdt}{
For arbitrary matrices 
$A\in\CC^{n,n}$, $B\in\CC^{n,m}$, 
$Q=Q'\in\CC^{n+m,n+m}$
the following conditions are equivalent:
\begin{itemize}
\item[(a)] there exists 
$P=P'\in\CC^{n,n}$
such that the Hermitian form $\s_P:\ \CC^n\times\CC^m\mapsto\RR$
defined by 
\EQ{k5}{  \s_P(x,u)=\s(x,u)+x'Px-(Ax+Bu)'P(Ax+Bu)}
is positive definite;
\item[(b)] the Hermitian form $\s$ is positive definite on the
subspace 
\EQ{k6}{  \cl L(z)=\{(x,u)\in\CC^n\times\CC^m:\ zx=Ax+Bu\}}
for all $z\in\TT$.
\end{itemize}
Moreover, when matrices $A,B,Q$ in (b) are real,
the corresponding matrix $P$ from (a) can be chosen to be
real as well.
}

\subsubsection{Non-Strict Linear Matrix Inequalities}

Since $\s_P(x,u)=\s(x,u)$ for $(x,u)\in \cl L(z)$, $z\in\TT$,
existence of a $P=P'\in\CC^{n,n}$ for which the Hermitian form
\rf{k5} is positive semidefinite implies that $\s$ is 
positive semidefinite on $\cl L(z)$ for all $z\in\TT$.
In general, the inverse implication is not
true: for example, when
$A=0$, $B=0$, and $\s(x,u)=\tx{Re}(x'u)$, the subspace $\cl L(z)$,
for all $z\in\TT$,
consists of all pairs $(0,u)$ with $u\in\CC$, and, accordingly,
$\s(x,u)=0$ for $(x,u)\in\cl L(z)$, $z\in\TT$. However, there exists
no $P=P'\in\CC^{1,1}$ (i.e. $P\in\RR$) for which
$\s_P(x,u)=\tx{Re}(x'u)+P|x|^2$ is positive semidefinite.

The following statement is a "non-strict LMI" 
version of the KYP Lemma which trades strict positivity for
controllability of the pair $(A,B)$. Recall that a pair $(A,B)$
of matrices $A\in\CC^{n,n}$, $B\in\CC^{n,m}$ 
is called {\sl controllable} when the matrix 
$[\lambda I_n-A,\ B]$ is right
invertible for all $\lambda\in\CC$.

\THM{thm:kyplminsdt}{
Assume that the pair $(A,B)$ of matrices
$A\in\CC^{n,n}$, $B\in\CC^{n,m}$ 
is controllable.  Then
for every matrix
$Q=Q'\in\CC^{n+m,n+m}$
the following conditions are equivalent:
\begin{itemize}
\item[(a)] there exists 
$P=P'\in\CC^{n,n}$
such that the Hermitian form $\s_P:\ \CC^n\times\CC^m\mapsto\RR$
defined by \rf{k5} is positive semidefinite;
\item[(b)] the Hermitian form $\s$ is positive semidefinite on the
subspace $\cl L(z)$ defined by
\rf{k6}
for all $z\in\TT$.
\end{itemize}
Moreover, when matrices $A,B,Q$ in (b) are real,
the corresponding matrix $P$ from (a) can be chosen to be
real as well.
}

\subsection{KYP Lemma in Continuous Time}
Continuous time (CT) versions of the KYP lemma are similar
to their DT counterparts.

\THM{thm:kypstabct}{
Assume that the pair $(A,B)$, where 
$A\in\CC^{n,n}$ and $B\in\CC^{n,m}$, 
is stabilizable, in the sense that
there exists a matrix $K\in\CC^{m,n}$ such that $sI_n-A-BK$ is
invertible for all $s\in\CC$ with $\tx{Re}(s)\ge0$.
Then for every matrix
$Q=Q'\in\CC^{n+m,n+m}$
the following conditions are equivalent:
\begin{itemize}
\item[(a)] there exist matrices
$P=P'\in\CC^{n,n}$, $C\in\CC^{m,n}$, $D\in\CC^{m,m}$ 
such that 
\EQ{k12}{\s(x,u)-2\tx{Re}[x'P(Ax+Bu)]=|Cx+Du|^2\ \ \ 
\forall\ x\in\CC^n,\ u\in\CC^m,}
\EQ{k13}{ 
\det\AR{cc}{A-sI_n & B\\ C & D}\neq0\ \ 
\forall\ s\in\CC_+;}
\item[(b)]  the matrix $\Pi(s)$ defined by \rf{k4}
is positive definite for all $s\in j\RR$
except, possibly,
a finite subset.
\end{itemize}
Moreover, when matrices $A,B,Q$ in (b) are real,
the corresponding matrices $P,C,D$ from (a) can be chosen to be
real as well.
}

\THM{thm:kyplmisct}{
For arbitrary matrices 
$A\in\CC^{n,n}$, $B\in\CC^{n,m}$, 
$Q=Q'\in\CC^{n+m,n+m}$
the following conditions are equivalent:
\begin{itemize}
\item[(a)] there exists 
$P=P'\in\CC^{n,n}$
such that the Hermitian form $\s_P:\ \CC^n\times\CC^m\mapsto\RR$
defined by 
\EQ{k15}{  \s_P(x,u)=\s(x,u)+2\tx{Re}[x'P(Ax+Bu)]}
is positive definite;
\item[(b)] the Hermitian form $\s$ is positive definite on the
subspace $\cl L(s)$ 
for all $s\in j\RR\cup\{\infty\}$, where $\cl L(z)$ is defined by \rf{k6} 
for $z\in\CC$, and
\EQ{k15x}{\cl L(\infty)=\{0\}\times\CC^m.}
\end{itemize}
Moreover, when matrices $A,B,Q$ in (b) are real,
the corresponding matrix $P$ from (a) can be chosen to be
real as well.
}

\THM{thm:kyplminsct}{
Assume that the pair $(A,B)$ of matrices
$A\in\CC^{n,n}$, $B\in\CC^{n,m}$ 
is controllable. Then
for every matrix
$Q=Q'\in\CC^{n+m,n+m}$
the following conditions are equivalent:
\begin{itemize}
\item[(a)] there exists 
$P=P'\in\CC^{n,n}$
such that the Hermitian form $\s_P:\ \CC^n\times\CC^m\mapsto\RR$
defined by \rf{k15} is positive semidefinite;
\item[(b)] the Hermitian form $\s$ is positive semidefinite on the
subspace $\cl L(s)$ defined by
\rf{k6}
for all $s\in j\RR$.
\end{itemize}
Moreover, when matrices $A,B,Q$ in (b) are real,
the corresponding matrix $P$ from (a) can be chosen to be
real as well.
}

\section{A Minimax Theorem}
It is easy to show that the inequality
\EQ{k21}{  \inf_{v\in V}\sup_{w\in W}g(v,w)\ge
\sup_{w\in W}\inf_{v\in V}g(u,w)}
holds for arbitrary sets $V,W$ and arbitrary
real-valued function $g:\ V\times W\mapsto\RR$.
The term {\sl minimax theorem} refers to a family of statements
providing conditions (usually involving
convexity of $g$ with respect to $v$ and
concavity of $g$ with respect to $w$)
under which the inequality in \rf{k21}
is actually an equality, i.e.
\EQ{k24}{  \inf_{v}\sup_{w}g(v,w)=
\sup_{w}\inf_{v}g(v,w).}

In this section, we are particularly interested in 
a specific minimax statement partially motivated by 
the KYP Lemma.

\subsection{Minimax Theorems for Discrete Time LTI Systems }
For a positive integer $m$ let $\ell^2_m$ denote the standard real
Hilbert space of all
one-sided real $m$-vector valued square summable sequences, 
i.e. functions $u:\ \ZZ_+\mapsto\RR^m$ such that
\EQ{k22}{  \|u(\cdot)\|^2=\sum_{t=0}^\infty|u(t)|^2<\infty.} 

Given a Schur matrix $A\in\RR^{n,n}$ (i.e. such that
$zI_n-A$ is not singular for $|z|\ge1$), a vector $a\in\RR^n$,
and matrices $B_1\in\RR^{n,k}$, $B_2\in\RR^{n,q}$,
$Q\in\RR^{n+k+q,n+k+q}$, consider the functional 
$g:\ \ell^2_k\times\ell^2_q\mapsto\RR$ defined by
\EQ{k23}{ g(v(\cdot),w(\cdot))=\sum_{t=0}^\infty\s(x(t),v(t),w(t)):\ \ 
x(t+1)=Ax(t)+B_1v(t)+B_2w(t),\ \ \ x(0)=a,}
where
\EQ{k25}{  \s(x,v,w)=\AR{c}{x\\ v\\ w}'Q\AR{c}{x\\ v\\ w}\ \ \ 
(x\in\RR^n,\ v\in\RR^k,\ w\in\RR^q).}
Consider also the associated matrix $\Pi=\Pi(z)$ defined by \rf{k4}
with $B=[B_1,\ B_2]$, and its partition
\EQ{k26}{  \Pi(z)=\AR{cc}{\Pi_{11}(z) & \Pi_{12}(z)\\ 
\Pi_{21}(z) & \Pi_{22}(z)},\ \ 
\Pi_{11}(z)\in\CC^{k,k},\ \Pi_{22}(z)\in\CC^{q,q}.}

Our objective is to formulate conditions, in terms of matrices
$\Pi_{ij}$, which guarantee that equality
\rf{k24} is satisfied for
all $a\in\RR^n$ for the functional 
$g:\ \ell^2_k\times\ell^2_q\mapsto\RR$ defined by \rf{k23}.
We are also interested in formulating conditions which ensure
that the associated partial infimum and supremum
\EQ{k28}{ g_v(w)=\inf_{v}g(v,w),\ \ \ 
g_w(v)=\sup_{w}g(v,w)}
are finite, and that
the resulting functions $g_v:\ \ell^2_q\mapsto\RR$, 
$g_w:\ \ell^2_k\mapsto\RR$ are continuous
in the standard Hilbert space metrics 
of $\ell^2_q$ and $\ell^2_k$, respectively.

\subsubsection{A Counterexample}
\label{sec:minmaxcount}

The Parceval identity can be used to show that $g(v,w)$ 
from \rf{k23} is convex
with respect to $v$ if and only if $\Pi_{11}(z)\ge0$ for all $z\in\TT$.
Similarly, $g(v,w)$ is concave
with respect to $w$ if and only if $\Pi_{22}(z)\le0$ for all $z\in\TT$.
However, these assumptions are far from being sufficient to assure that
the
minimax identity \rf{k24} is satisfied, 
as demonstrated by the example with
\[  A=\AR{cc}{0&0\\ 0&0},\ B_1=\AR{c}{1\\ 0},\ B_2=\AR{c}{0\\ 1},\ 
a=\AR{c}{1\\ 0},\ 
Q=\AR{cccc}{1&0&-1&0\\ 0&-1&1&1\\ -1&1&1&-1\\ 0&1&-1&-1},\]
i.e. when 
\[ g(v(\cdot),w(\cdot))=
\sum_{t=0}^\infty\{|v(t)-x_1(t)|^2-2v(t)[w(t)-x_2(t)]-
|w(t)-x_2(t)|^2\},\] 
\[ \tx{subject to}\ \ 
x_1(t+1)=v(t),\ x_2(t+1)=w(t),\ x_1(0)=1,\ x_2(0)=0,\]
and
\[  \Pi(z)=\AR{cc}{|z-1|^2 & z-1\\ z'-1 & -|z-1|^2}.\]

Using the fact that the set of all possible sequences $w-x_2$ is dense
in $V=W=\ell^2$, we conclude that
\[  \sup_{w\in W}g(v,w)=\sum_{t=0}^\infty\{|v(t)-x_1(t)|^2+|v(t)|^2\}
\ge|v(0)-1|^2+|v(0)|^2\ge0.5,\]
and hence the left side in \rf{k24} is not smaller than 0.5.

On the other hand, using the fact that
\[  \sum_{t=0}^\infty\{x_1(t)x_2(t)-v(t)w(t)\}=
\sum_{t=0}^\infty\{x_1(t)x_2(t)-x_1(t+1)x_2(t+1)\}=x_1(0)x_2(0)=0,\]
one can re-write the sum for $g$ as
\[ g(v(\cdot),w(\cdot))=
\sum_{t=0}^\infty\{|v(t)-x_1(t)|^2+
2[v(t)-x_1(t)]x_2(t)-|w(t)-x_2(t)|^2\}.\]
Since the set of all possible sequences $v-x_1$ is dense
in $\ell^2$, we conclude that
\[  \inf_{v\in V}g(v,w)=
\sum_{t=0}^\infty\{-|x_2(t)|^2-|w(t)-x_2(t)|^2\},\]
and hence the right side in \rf{k24} is zero: the 
minimax equality does not
hold in this case.

It is instructive to note that, in this case, both functions
$g_v$ and $g_w$ from \rf{k28}
are finite and continuous in the standard Hilbert space
topology of $\ell^2$,
but the minimax identity is still not valid.

\subsubsection{A Sufficient Condition of Minimax}
The following statement shows that the minimax identity
\rf{k24} holds for the functional \rf{k23} when there exist
$\e>0$ and $z_0\in\TT$ such that
\EQ{k27}{  \AR{cc}{\Pi_{11}(z) & \e\Pi_{12}(z)\\ 
\e\Pi_{21}(z) & -\Pi_{22}(z)}\ge0\ \  \ \forall\ z\in\TT,\  \ 
\Pi_{11}(z_0)>0,\ \Pi_{22}(z_0)<0.}

\THM{thm:minmaxdt}{
Let  $A\in\RR^{n,n}$ be a Schur matrix. Assume that
matrices $B_1\in\RR^{n,k}$, $B_2\in\RR^{n,q}$,
$Q\in\RR^{n+k+q,n+k+q}$ are such that condition \rf{k27},
where $\Pi_{ij}$ are defined by \rf{k26} 
and \rf{k4} with $B=[B_1,B_2]$,
is satisfied for some $\e>0$ and $z_0\in\TT$.
Then 
for every $a\in\RR^n$ for the functional 
$g:\ \ell^2_k\times\ell^2_q\mapsto\RR$ defined by
\rf{k23},\rf{k25}
\begin{itemize}
\item[(a)] the partial optimal values
in \rf{k28} are continuous in the standard norm topologies
of $\ell^2_q$ and $\ell^2_k$;
\item[(b)] the minimax identity in \rf{k24} 
is satisfied.
\end{itemize}
}

\subsubsection{Minimax and Integral Quadratic Constraints}

For a positive integer $m$ let $\ell_m$ denote the set
of all functions $u:\ \ZZ_+\mapsto\RR^m$ (in particular,
$\ell_m^2$ is a subset of $\ell_m$).
In modeling discrete time dynamical systems, $m$-dimensional
{\sl signals} can be represented by the elements of $\ell_m$.
Accordingly, a DT {\sl system}
$\D$ with $k$-dimensional input $v$ and $q$-dimensional
output $w$ is viewed as a subset $\D\subset\ell_k\times\ell_q$.
Let us call such system $\D\subset\ell_k\times\ell_q$
 {\sl weakly causally stable} if for 
every $T\in\ZZ_+$, $(v,w)\in\D$, and
$v_*\in\ell_k^2$ such that $v(t)=v_*(t)$ for all $t\le T$
there exists a sequence of
elements $(v_i,w_i)\in\D\cap(\ell^2_k\times\ell^2_q)$, 
such that $v_i(t)=v(t)$ and
$w_i(t)=w(t)$ for all
$t\le T$, and $\|v_i-v_*\|\to0$ as $i\to\infty$.

Given real matrices $A\in\RR^{n,n}$, 
$B_1\in\RR^{n,k}$, $B_2\in\RR^{n,q}$, $Q\in\RR^{n+k+q,n+k+q}$,
where $A$ is
a Schur matrix, and a subset $X_0\subset\RR^n$,
let us say that system $\D\subset\ell_k\times\ell_q$
satisfies the {\sl conditional} Integral Quadratic Constraint (IQC)
defined by $A$, $B_1$, $B_2$, $Q$, $X_0$ if there
exists a continuous function 
$\kappa:\ \RR^n\times\RR^k\times\RR^q\mapsto\RR$
such that
\EQ{k31}{\sum_{t=0}^\infty \s(x(t),v(t),w(t))\ge-\kappa(x_0,v(0),w(0))} 
for all 
$(v,w)\in\D\cap(\ell^2_k\times\ell^2_q)$, $x_0\in X_0$,
where $x(\cdot)$ is defined by
$v(\cdot)$, $w(\cdot)$, and $x_0$ according to
\EQ{k32}{x(t+1)=Ax(t)+B_1v(t)+b_2w(t),\ \ \ x(0)=x_0.}
Similarly, 
let us say that $\D\subset\ell_k\times\ell_q$
satisfies the {\sl complete} IQC
defined by $A$, $B_1$, $B_2$, $Q$, $X_0$ if there
exists a continuous function 
$\kappa:\ \RR^n\times\RR^k\times\RR^q\mapsto\RR$
such that
\EQ{k33}{\sum_{t=0}^T \s(x(t),v(t),w(t))\ge-\kappa(x_0,v(0),w(0))\ \ \ 
(T\ge0)} 
for all  
$(v,w)\in\D$, $x_0\in X_0$, and 
$x\in\ell_n$ satisfying \rf{k32}.

An important step in the IQC framework of nonlinear system analysis
is to establish that a particular conditional IQC \rf{k31}
implies the corresponding complete IQC \rf{k33}. The
implication is not always
true: for example, when 
\[   \D=\{(v,w)\in\ell\times\ell:\ \ \ 
w(t+1)=v(t)\ \forall\ t\in\ZZ_+\},\]
\[  A=B_1=B_2=0,\ \ X_0=\{0\},\ \ \s(x,v,w)=|w|^2-|v|^2\]
then the conditional IQC \rf{k31} is satisfied with
$\kappa(x_0,v_0,w_0)=|v_0|^2$, but the associated complete IQC 
\rf{k33} does not take place for any function $\kappa$.

The following statement, based on the
minimax identity established in Theorem~\ref{thm:minmaxdt},
 provides sufficient conditions,
expressed in terms of matrices $A$, $B$, and $Q$,
under which the conditional IQC from \rf{k31} implies the complete
IQC from \rf{k33}.

\THM{thm:minmaxiqcdt}{
Let $\D\subset\ell_k\times\ell_q$
be a  weakly causally stable system which satisfies the conditional
IQC defined by real matrices $A\in\RR^{n,n}$, 
$B_1\in\RR^{n,k}$, $B_2\in\RR^{n,q}$, $Q\in\RR^{n+k+q,n+k+q}$,
where $A$ is
a Schur matrix, and a subset $X_0\subset\RR^n$.
Assume that
\begin{itemize}
\item[(a)] condition \rf{k27},
where $\Pi_{ij}$ are defined by \rf{k26} 
and \rf{k4} with $B=[B_1,B_2]$,
is satisfied for some $\e>0$ and $z_0\in\TT$;
\item[(b)]
there exist real matrices $C\in\RR^{k,n}$, $D_1\in\RR^{k,k}$,
$D_2\in\RR^{k,q}$ such that the quadratic form $\s$ defined by \rf{k25}
satisfies the inequality
\EQ{k34}{\s(x,v,w)\le|Cx+D_1v+D_2w|^2\ \ \ 
\forall\ x\in\RR^n,\ v\in\RR^k,\ w\in\RR^q,}
and 
\EQ{k35}{  \det\AR{cc}{\lambda A-I & \lambda B_1\\ C & D_1}\neq0\ \ \ 
\forall\ \lambda\in\CC,\ |\lambda|<1.}
\end{itemize}
Then $\D$ satisfies the complete
IQC defined by $A$,$B_1$, $B_2$, $Q$, $X_0$.
}

\subsection{Minimax Theorems for Continuous Time LTI Systems }
For a positive integer $m$ let $L^2_m$ denote the standard real
Hilbert space of all
 real $m$-vector valued square integrable functions
$u:\ [0,\infty)\mapsto\RR^m$, equipped with the norm
\EQ{k22x}{  \|u(\cdot)\|^2=\int_{0}^\infty|u(t)|^2dt<\infty.} 

Given a Hurwitz matrix $A\in\RR^{n,n}$ (i.e. such that
$sI_n-A$ is not singular for $\tx{Re}(s)\ge0$), a vector $a\in\RR^n$,
and matrices $B_1\in\RR^{n,k}$, $B_2\in\RR^{n,q}$,
$Q\in\RR^{n+k+q,n+k+q}$, consider the functional 
$g:\ L^2_k\times L^2_q\mapsto\RR$ defined by
\EQ{k23x}{ g(v(\cdot),w(\cdot))=\int_{0}^\infty\s(x(t),v(t),w(t))dt:\ \ 
\dot x(t)=Ax(t)+B_1v(t)+B_2w(t),\ \ \ x(0)=a,}
where $\s(\cdot)$ is defined by \rf{k25}.
Consider also the associated matrix $\Pi(\cdot)$ defined by \rf{k4}
with $B=[B_1,\ B_2]$, and its partition \rf{k26}.

Our objective is to formulate conditions, in terms of matrices
$\Pi_{ij}$, which guarantee that equality
\rf{k24} is satisfied for
all $a\in\RR^n$ for the functional 
$g:\ L^2_k\times L^2_q\mapsto\RR$ defined by \rf{k23x}.
We are also interested in formulating conditions which ensure
that the associated partial infimum and supremum \rf{k28}
are finite, and that
the resulting functions $g_v:\ L^2_q\mapsto\RR$, 
$g_w:\ L^2_k\mapsto\RR$ are continuous
in the standard Hilbert space metrics 
of $L^2_q$ and $L^2_k$, respectively.

\subsubsection{A Sufficient Condition of Minimax}
The following statement shows that the minimax identity
\rf{k24} holds for the functional \rf{k23x} when there exist
$\e>0$ and $s_0\in j\RR$ such that
\EQ{k27x}{  \AR{cc}{\Pi_{11}(s) & \e\Pi_{12}(s)\\ 
\e\Pi_{21}(s) & -\Pi_{22}(s)}\ge0\ \  \ \forall\ s\in j\RR,\  \ 
\Pi_{11}(s_0)>0,\ \Pi_{22}(s_0)<0.}

\THM{thm:minmaxct}{
Let  $A\in\RR^{n,n}$ be a Hurwitz matrix. Assume that
matrices $B_1\in\RR^{n,k}$, $B_2\in\RR^{n,q}$,
$Q\in\RR^{n+k+q,n+k+q}$ are such that condition \rf{k27x},
where $\Pi_{ij}$ are defined by \rf{k26} 
and \rf{k4} with $B=[B_1,B_2]$,
is satisfied for some $\e>0$ and $s_0\in j\RR$.
Then 
for every $a\in\RR^n$ for the functional 
$g:\ L^2_k\times L^2_q\mapsto\RR$ defined by
\rf{k23x},\rf{k25}
\begin{itemize}
\item[(a)] the partial optimal values
in \rf{k28} are continuous in the standard norm topologies
of $L^2_q$ and $L^2_k$;
\item[(b)] the minimax identity in \rf{k24} 
is satisfied.
\end{itemize}
}

\subsubsection{Minimax and Continuous Time  IQC}

For a positive integer $m$ let $L_m$ denote the set
of all locally square integrable functions 
$u:\ [0,\infty)\mapsto\RR^m$ (in particular,
$L_m^2$ is a subset of $L_m$).
In modeling continuous time dynamical systems, $m$-dimensional
{\sl signals} can be represented by the elements of $L_m$.
Accordingly, a CT {\sl system}
$\D$ with $k$-dimensional input $v$ and $q$-dimensional
output $w$ is viewed as a subset $\D\subset L_k\times L_q$.
Let us call such system $\D\subset L_k\times L_q$
 {\sl weakly causally stable} if for 
every $T\ge0$, $(v,w)\in\D$, and
$v_*\in L_k^2$ such that $v(t)=v_*(t)$ for all $t\le T$
there exists a sequence of
elements $(v_i,w_i)\in\D\cap(L^2_k\times L^2_q)$, 
such that $v_i(t)=v(t)$ and
$w_i(t)=w(t)$ for all
$t\le T$, and $\|v_i-v_*\|\to0$ as $i\to\infty$.

Given real matrices $A\in\RR^{n,n}$, 
$B_1\in\RR^{n,k}$, $B_2\in\RR^{n,q}$, $Q\in\RR^{n+k+q,n+k+q}$,
where $A$ is
a Hurwitz matrix, and a subset $X_0\subset\RR^n$,
let us say that system $\D\subset\ell_k\times\ell_q$
satisfies the {\sl conditional} Integral Quadratic Constraint (IQC)
defined by $A$, $B_1$, $B_2$, $Q$, $X_0$ if there
exists a continuous function 
$\kappa:\ \RR^n\times\RR^k\times\RR^q\mapsto\RR$
such that
\EQ{k31x}{\int_{0}^\infty \s(x(t),v(t),w(t))dt\ge-\kappa(x_0,v(0),w(0))} 
for all 
$(v,w)\in\D\cap(L^2_k\times L^2_q)$, $x_0\in X_0$,
where $x(\cdot)$ is defined by
$v(\cdot)$, $w(\cdot)$, and $x_0$ according to
\EQ{k32x}{\dot x(t)=Ax(t)+B_1v(t)+B_2w(t),\ \ x(0)=x_0,}
which is understood, in a generalized sense, as
\[x(t)=x_0+\int_0^t[Ax(\tau)+B_1v(\tau)+B_2w(\tau)]d\tau\ (t\ge0).\]
Similarly, 
let us say that $\D\subset L_k\times L_q$
satisfies the {\sl complete} IQC
defined by $A$, $B_1$, $B_2$, $Q$, $X_0$ if there
exists a continuous function 
$\kappa:\ \RR^n\times\RR^k\times\RR^q\mapsto\RR$
such that
\EQ{k33x}{\int_{0}^T \s(x(t),v(t),w(t))dt\ge-\kappa(x_0,v(0),w(0))\ \ \ 
(T\ge0)} 
for all  
$(v,w)\in\D$, $x_0\in X_0$, and 
$x\in L_n$ satisfying \rf{k32x}.

\THM{thm:minmaxiqcct}{
Let $\D\subset L_k\times L_q$
be a  weakly causally stable system which satisfies the conditional
IQC defined by real matrices $A\in\RR^{n,n}$, 
$B_1\in\RR^{n,k}$, $B_2\in\RR^{n,q}$, $Q\in\RR^{n+k+q,n+k+q}$,
where $A$ is
a Hurwitz matrix, and a subset $X_0\subset\RR^n$.
Assume that
\begin{itemize}
\item[(a)] condition \rf{k27x},
where $\Pi_{ij}$ are defined by \rf{k26} 
and \rf{k4} with $B=[B_1,B_2]$,
is satisfied for some $\e>0$ and $s_0\in j\RR$;
\item[(b)]
there exist real matrices $C\in\RR^{k,n}$, $D_1\in\RR^{k,k}$,
$D_2\in\RR^{k,q}$ such that the quadratic form $\s$ defined by \rf{k25}
satisfies the inequality \rf{k35},
and 
\EQ{k35x}{  \det\AR{cc}{A-sI & B_1\\ C & D_1}\neq0\ \ \ 
\forall\ s\in\CC_+.}
\end{itemize}
Then $\D$ satisfies the complete
IQC defined by $A$,$B_1$, $B_2$, $Q$, $X_0$.
}

\section{Appendix}
This section contains proofs of main statements made in the paper,
including a brief description of some classical mathematical constructions used in the proofs.

\subsection{Quadratic Optimization and Minimax}
We begin by summarizing some elementary statements concerning
quadratic functionals and real Hilbert spaces.

\subsubsection{Quadratic Forms}

A function $\s:\ V\mapsto\RR$ defined on a real vector space
$V$ is called a {\sl quadratic form}
when $\s(v)=b(v,v)$ for all $v\in V$, where
$b:\ V\times V\mapsto\RR$ is a {\sl symmetric bilinear} function, i.e.
\EQ{k025}{ b(u,v)=b(v,u),\ \ b(u,xv+yw)=xb(u,v)+yb(u,w)\ \ 
\forall\ u,v,w\in V,\ x,y\in\RR.}
This correspondence between symmetric bilinear functions and
quadratic forms is a bijection, as 
$b(\cdot,\cdot)$ can be recovered from $\s(\cdot)$ according to the
identity
\[ b(u,v)=\frac{\s(u+v)-\s(u-v)}4.\]
The quadratic form $\s$ is called {\sl positive definite} 
(notation $\s\gg0$) when $\s(v)>0$ for all
$v\neq0$, and {\sl positive semidefinite} 
(notation $\s\ge0$) when $\s(v)\ge0$ for all $v\in V$.
Due to the identity
\EQ{k51}{ t\s(v)+(1-t)\s(u)-\s(tv+(1-t)u)=t(1-t)\s(v-u)\ \ 
\forall\ v,u\in V,\ t\in\RR,}
which is valid for every quadratic form $\s:\ V\mapsto\RR$,
$\s$ is convex if and only if it is positive semidefinite.

For example, a symmetric real matrix $Q=Q'\in\RR^{n,n}$
defines a symmetric bilinear form $b_Q:\ \RR^n\times\RR^n\mapsto\RR$
according to $b_Q(v,u)=v'Qu$, and the associated quadratic form
$\s_Q(v)=b_Q(v,v)$; the form $\s_Q$ (equivalently, the
matrix $Q=Q'$) is  positive definite (or semidefinite)
when all eigenvalues of $Q$ are positive (notation $Q>0$)
or, respectively, non-negative (notation $Q\ge0$).

\subsubsection{Quadratic Optimization
and Real Hilbert Spaces}

In this paper, the terminology of quadratic forms is used 
to {\sl formulate} statements (this makes assumptions
easier to verify in applications), while the more flexible
Hilbert space viewpoint
is employed in the corresponding {\sl proofs}. The definitions and
statements of this subsection
facilitate easy switching between the two frameworks.

Let $b:\ V\times V\mapsto\RR$ be a symmetric bilinear form 
on a real vector space
$V$ such that the corresponding quadratic form 
$\s(v)=b(v,v)$ is positive definite.  Since
the quadratic function 
\[ t\in\RR\ \mapsto\ \s(v+tu)=\s(v)+2tb(v,u)+t^2\s(u)\]
takes only non-negative values, its discriminant is not positive, which
yields the Cauchy-Schwartz inequality
\EQ{CAL2}{  |b(v,u)|^2\le\s(v)\s(u)\ \ \ \forall\ v,u\in V,}
and in turn implies that
the function $v\mapsto|v|_{\s}=\s(v)^{1/2}$
is a {\sl norm} on $V$. 

Let $V^{\s}$ be  the set of all linear
functions $f:\ V\mapsto\RR$ such that
\[  |f|_{\s}\bydef\sup\{f(v):\  \s(v)\le1\}<\infty.\]
As a dual of a normed space $(V,|\cdot|_{\s})$, the pair
$(V^{\s},|\cdot|_{\s})$ defines a Banach space. 
Let  $\pi_{\s}:\ V\mapsto V^{\s}$ be the
"natural correspondence" mapping
every $v\in V$ to $f=\pi_{\s}v\in V^{\s}$ according to 
$f(u)=b(v,u)$. 

The quantity $|f|_{\s}^2$, where $f\in V^{\sharp}$ 
can also be interpreted
as the minimal upper bound in the
{\sl quadratic optimization} task
\EQ{CAL3}{ 2fv-\s(v)\mapsto\sup_{v\in V},}
because
\[ \inf_v\{2f(v)-\s(v)\}=\inf_{\s(v)\le1}\inf_{t\in\RR}\{2f(tv)-\s(tv)\}
=\inf_{\s(v)\le1}\inf_{t\in\RR}\{2f(v)t-\s(v)t^2\}.\]

\THM{thm:switch}{
Let $b:\ V\times V\mapsto\RR$ be a symmetric bilinear form 
on a real vector space
$V$ such that the corresponding quadratic form $\s(v)=b(v,v)$ is positive definite. Then
\begin{itemize}
\item[(a)] the set $\pi_{\s}V$ is dense in $(V^{\s},|\cdot|_{\s})$;
\item[(b)] there exists a (unique) symmetric bilinear form 
$\bar b:\ V^{\s}\times V^{\s}\mapsto\RR$ such that
$|f|_{\s}^2=\bar b(f,f)$ for all $f\in V^{\s}$;
\item[(c)] for the bilinear form $\bar b$ defined in (b),
the identity $f(v)=\bar b(f,\pi_{\s}v)$ holds for all $f\in V^{\s}$ and
$v\in V$.
\end{itemize}
}

Theorem~\ref{thm:switch} establishes $(V^{\s},|\cdot|_{\s})$
as a (real) Hilbert space, and provides a linear 
norm-preserving bijection
$\pi_{\s}$
between vectors from $V$
and elements of a dense subspace $\pi_{\s}V$ of $V^{\s}$.
It also shows that the minimal upper bound in
  quadratic optimization \rf{CAL3}, as a function of $f\in V^\sharp$,
 is a
quadratic form on the subset  $V^{\s}$ where its values are finite.

In this paper, we will use either $|w|$ or $|w|_H$
to denote the norm of a vector $w$ in a Hilbert space $H$.
In addition, the shortcut
$v'u$ will denote the scalar product of
two vectors $v,u$ from the same Hilbert space $H$.
This notation can be motivated by the natural
association of vectors $v\in H$ 
with bounded linear functions $L_v:\ \RR\mapsto H$
defined by $L_v(t)=tv$. Accordingly, the adjoint $v'$ 
is the linear function $v':\ H\mapsto\RR$ mapping
$u$ to the scalar product of $v$ and $u$, and the composition
$v'u$, where $v,u\in H$, is a linear function mapping $\RR$ to $\RR$,
i.e. a real number, which equals the scalar product of $v$ and $u$

\subsubsection{Quadratic Minimax}

The following statement 
lists sufficient conditions for the minimax identity in quadratic
optimization.

\THM{thm:qminmax}{
Let $V,W$ be real vector spaces. Let $g:\ V\times W\mapsto\RR$
be defined by
\EQ{CAL21}{  g(v,w)=\s(v)+2p(v,w)-\mu(w)-2f(v)+2h(w)+r,}
where
$\s:\ V\mapsto\RR$, 
$\mu:\ W\mapsto\RR$, $p:\ V\times W\mapsto\RR$,
$f:\ V\mapsto\RR$, 
$h:\ W\mapsto\RR$, and $r\in\RR$
are two positive definite quadratic forms,
a bilinear functional, two linear functions, and a real number. 
Assume that
\begin{itemize}
\item[(i)] there exists $c\ge0$ such that
$c^2\s(v)\mu(w)\ge |p(v,w)|^2$ for all $v\in V$, $w\in W$;
\item[(ii)] $\inf_{v\in V}g(v,0)>-\infty$ and 
$\sup_{w\in W}g(0,w)<+\infty$.
\end{itemize}
Then 
\begin{itemize}
\item[(a)] the minimax equality \rf{k24} holds;
\item[(b)] there exists a constant $c_1\ge0$ such that
\[ \inf_{v\in V} g(v,w)\ge-c_1(1+\mu(w))\ \ \forall w\in W,\ \ \ 
\sup_w g(v,w)\le c_1(1+\s(v))\ \ \forall\ v\in V.\]
\end{itemize}
}

\begin{pf}
By (i), for every $v\in V$ the 
function $f_v:\ W\mapsto\RR$ defined by
$f_v(w)=p(v,w)$ is linear and satisfies 
$|f_v|_{\mu}\le c|v|_{\s}$, i.e. $f_v\in W^{\mu}$.
Since the corresponding function
$\pi_{\s}v\in V^{\s}\mapsto f_v\in W^\mu$
is linear and bounded, it can be extended to a bounded linear
operator $L:\ V^{\s}\mapsto W^\mu$ such that
\[  p(v,w)=\bar w'L\bar v\ \ 
(\bar w=\pi_\mu w,\bar v=\pi_{\s}v)\ \ \forall\ v\in V,w\in W.\]
Since $L$ is bounded, its adjoint $L'$ is well defined and bounded
as well.
Also, by (ii), $f\in V^{\s}$ and $h\in W^\mu$, hence the identity
$g(v,w)=\bar g(\bar v,\bar w)$ holds for
$\bar w=\pi_\mu w$, $\bar v=\pi_{\s}v$, and
\[ \bar g(\bar v,\bar w)=
|\bar v|^2+2\bar w'L\bar v-|\bar w|^2-2f'\bar v+2h'\bar w+r.\]

Let $A:\ V^{\s}\times W^{\mu}\mapsto V^{\s}\times W^{\mu}$
be the linear operator with block representation
\[  A=\AR{cc}{I & L'\\ -L & I},\]
i.e. $A(v,w)=(v+L'w,w-Lv)$. Since $A$ is bounded and $A+A'=2I$
is strictly positive definite,
$A$ must be invertible, and hence there exist $v_0\in V^{\s}$,
$w_0\in W^{\mu}$ such that $A(v_0,w_0)=(f,h)$.

Since $L$ is bounded and the subsets
$\pi_{\s}V$, $\pi_{\mu}W$ are dense in $V^{\s}$ and $W^{\mu}$
respectively, we have (using notation $\bar v=\pi_{\s}v$,
$\bar w=\pi_\mu(w)$ for $v\in V$ and $w\in W$):
\begin{eqnarray*}
\sup_{w\in W}g(v,w)&=&
  \sup_{\bar w\in \pi_{\mu}W}\bar g(\bar v,\bar w)\\
&=&\sup_{\hat w\in W^{\mu}}\bar g(\bar v,\hat w)\\
&=&|\bar v|^2+|L\bar v+h|^2-2f'\bar v+r\\
&=&|\bar v|^2+|L\bar v-Lv_0+w_0|^2-2(v_0+L'w_0)'\bar v+r\\
&=&|\bar v-v_0|^2+|L(\bar v-v_0)|^2+|w_0|^2-|v_0|^2-2w_0'Lv_0+r,
\end{eqnarray*}
hence the second inequality in (b) holds, and
\begin{eqnarray*}
\inf_{v\in V}\sup_{w\in W}g(v,w)&=& \inf_{\bar v\in \pi_{\s}V}
\sup_{\bar w\in \pi_{\mu}W}\bar g(\bar v,\bar w)\\
&=&|w_0|^2-|v_0|^2-2w_0'Lv_0+r.
\end{eqnarray*}

Similarly,
\begin{eqnarray*}
\inf_{v\in V}g(v,w)&=&
  \inf_{\bar v\in \pi_{\s}V}\bar g(\bar v,\bar w)\\
&=&\inf_{\hat v\in V^{\s}}\bar g(\hat v,\bar w)\\
&=&-|\bar w|^2-|L'\bar w-f|^2+2h'\bar w+r\\
&=&-|\bar w|^2-|L'\bar w-L'w_0-v_0|^2+2\bar w'(w_0-Lv_0) +r\\
&=&
-|\bar w-w_0|^2-|L'(\bar w-w_0)|^2+|w_0|^2-|v_0|^2-2w_0'Lv_0+r,
\end{eqnarray*}
hence the first inequality in (b) holds, and
\begin{eqnarray*}
\sup_{w\in W}\inf_{v\in V}g(v,w)&=&
\sup_{\bar w\in \pi_{\mu}W}\inf_{\bar v\in \pi_{\s}V}
\bar g(\bar v,\bar w)
=|w_0|^2-|v_0|^2-2w_0'Lv_0+r,
\end{eqnarray*}
which establishes the minimax identity.

The bounds from (b) follow from the explicit expressions for the
partial optimal values, and from the boundedness of $L$ and $L'$.

\end{pf}

\subsection{KYP Lemma Proofs}
This section contains proofs of the statements associated with the
KYP Lemma.

\subsubsection{Theorem~\ref{thm:kypstabdt}, (a)$\Rightarrow$(b)}
For $z\not\in\Lambda(A)$
let 
\[ L(z)=(zI-A)^{-1}B,\ \ \ \ H(z)=D+CL(z).\]
Substituting
$x=L(z)u$ (which means $Ax+Bu=zx$) with $z\in\TT$ into \rf{k2}
 yields
\EQ{CAL6}{ \Pi(z)=H(z)'H(z)\ \ \ 
\forall\ z\in\TT\backslash\Lambda(A),}
hence $\Pi(z)\ge0$ for 
$z\in\TT\backslash\Lambda(A)$. Moreover, since
\[  \det\AR{cc}{z^{-1}A-I & z^{-1}B\\ C & D}=\det H(z)\]
for $z\neq0$, $z\not\in\Lambda(A)$,
the rational function $z\mapsto\det H(z)$ is not identically
equal to zero, and hence $\det H(z)\neq0$ for all $z\in\CC$ except,
possibly, a finite subset. Hence \rf{CAL6} implies that
$\Pi(z)$ is positive definite for all $z\in\CC$ except,
possibly, a finite subset.

\subsubsection{Theorem~\ref{thm:kypstabdt}, (b)$\Rightarrow$(a)}
To prove the implication, we consider the associated
optimization setup \rf{k7}-\rf{k7B}, which can be recognized as
a case of quadratic optimization. The key step is to show that
the infimum in \rf{k7}-\rf{k7B} is finite. Then, according to
Theorem~\ref{thm:switch}, $\inf\Phi$ is a quadratic form
of $a$. We define $P=P'$ by $\inf\Phi=-a'Pa$, and use
the Bellman equation from dynamic programming to
show that conditions \rf{k2},\rf{k3} are satisfied.

\begin{itemize}
\item[(a)] Let $\ell_m^2$ be the set
of complex
square summable sequences $w:\ \ZZ_+\mapsto\CC$,
equipped with the natural structure of
a {\sl real} vector space. Since 
$A+BK$ is a Schur matrix, 
there is a linear one-to-one correspondence between the pairs
$(x,u)$ in \rf{k7A},\rf{k7B} and the pairs
$(w,a)\in\ell_m^2\times\CC^n$ which maps $(x,u)$ to
$(u-Kx,x(0))$. 

Using the Parceval identity, 
the functional $\Phi$ in \rf{k7} can be re-written in the form
\EQ{RIC7}{  \Phi=\int_{\TT}
\{ \hat w(z)'\a(z)\hat w(z)+2\tx{Re}[\hat w(z)\b(z)a]+a'\g(z)a\}dm(z),}
where
\[  \int_{\TT}f(z)dm(z)=\frac1{2\pi}\int_{-\pi}^{\pi}f(e^{j\t})d\t\]
denotes the standard Lebesque measure integral on the unit circle
$\TT$,
\[  \hat w(z)=\sum_{t=0}^\infty w(t)z^{-t}\]
is the Fourier transform of $w\in\ell_m^2$, 
a square integrable function $\hat w:\ \TT\mapsto\CC$, and
$\a,\b,\g$ are the rational matrix-valued
functions defined by the block decomposition
identity (to be satisfied for
$z\in\TT$)
\[  \AR{cc}{\a(z) & \b(z)'\\ \b(z) & \g(z)} = 
M(z)'QM(z),\]
with
\[
M(z)=\AR{cc}{I & 0\\ K & I}
\AR{cc}{(zI-A-BK)^{-1}&0\\ 0 & I}
\AR{cc}{B&I\\ I & 0}.\]
Since $A+BK$ is a Schur matrix, $\a,\b,\g$ have no poles on the
unit circle $\TT$. Also, since
\[ \a(z)=F(z)'\Pi(z)F(z),\ \ \tx{where}\ \ 
F(z)=[I-K(z-A)^{-1}B]^{-1}\]
for $z\in\TT$, the matrix $\a(z)$ is positive definite for all
$z\in\TT$ except, possibly, a finite subset, where it is positive
semidefinite. 

Since, at the points where $\a(z)$ is positive definite,
\[ \bar w'\a(z)\bar w+2\tx{Re}{\bar w'\b(z)a}\ge
- a'\b(z)'\a(z)^{-1}\b(z)a\ \ \ \forall\ \bar w\in\CC,\]
the conclusion $\inf\Phi>-\infty$ can be reached easily when
there exists a constant $c\in\RR$ such that
$\b(z)'\a(z)^{-1}\b(z)\le cI_m$ for all $z\in\TT$ with $\a(z)>0$.
While such $c\in\RR$ does not always exist, we can use the fact
that 
\[  \int_{z\in\TT} \hat w(z)'\d(z)dm(z)=0\]
for every $w\in\ell_m^2$ and every strictly proper rational
matrix $\d=\d(z)$ with no poles outside the open unit disk
$|z|<1$.

Indeed, to prove that $\inf\Phi>-\infty$, it is sufficient to
find a strictly proper rational matrix function $\d=\d(z)$
with no poles outside the open unit circle $|z|<1$, 
with the property that
there exists a constant $c\in\RR$ such that
\EQ{CAL8}{
(\b(z)-\d(z))'\a(z)^{-1}(\b(z)-\d(z))\le cI_m\ \ \ 
\tx{for}\  z\in\TT:\ \a(z)>0.}
Let
\[  R=\max_{z\in\TT}\lambda_{\max}(\a(z))\]
be the maximal eigenvalue of $\a(z)$ over $z\in\TT$
(it exists since $\a$ is continuous on $\TT$).
Then $\a(z)\ge\rho(z)I_m$ for all $z\in\TT$, where the scalar
rational function $\rho=\rho(z)$ is defined by
\[ \rho(z)=\det(\a(z))R^{1-m}.\]
Hence condition \rf{CAL8} will be satisfied, for some $c\in\RR$,
when the ratio $(\b-\d)/\rho$ is bounded on $\TT$,
i.e. when the unit circle zeroes of the scalar components
of $\b-\d$  match (counting multiplicity) the unit circle
zeroes of $\rho$. 

Recall that for every set of distinct complex numbers 
$(\lambda_i)_{i=1}^N$
and polynomials 
\[ p_i(\lambda)=\sum_{l=0}^{m_i-1}p_{i,l}\lambda^l\]
there exists a polynomial $p=p(\lambda)$ of degree
$\sum m_i$ such that 
\[ p(\lambda)-p_i(\lambda)=O((\lambda-\lambda_i)^{m_i})\ \ 
\tx{as}\ \ \lambda\to\lambda_i\ \  \ \forall\ i.\]
Hence the boundedness of $(\b-\d)/\rho$ on $\TT$ can be
achieved
by selecting $\d=\d(z)$ as a linear combination of a sufficiently
large number of monomials $z^{-i}$ with positive integer $i$,
which completes the proof of the inequality $\inf\Phi>-\infty$.

\item[(b)] 
Since $V(a)\bydef\inf\Phi>-\infty$ for every $a\in\CC^n$,
Theorem~\ref{thm:switch}, together with representation
\rf{RIC7}, imply that $V=V(a)$ is a quadratic form
of $a\in\CC^n$. Moreover, since multiplying a solution
$(x,u)$ of \rf{k7A} with $x(0)=a$ by $j$ yields a solution
$(jx,ju)$ of \rf{k7A} with $x(0)=ja$ and the same value
of $\Phi$, we have $V(ja)=V(a)$
for every $a\in\CC^n$, which implies that
$V(a)=-a'Pa$ for some fixed
complex $n$-by-$n$ matrix $P=P'$. The Bellman inequality
for the optimization task \rf{k7}-\rf{k7B} can be written in the form
\EQ{CAL9}{\inf_{u\in\CC^m}\{\s(x,u)+V(Ax+Bu)-V(x)\}=0.}
Since $\mu(x,u)\bydef\s(x,u)+V(Ax+Bu)-V(x)$ 
is a quadratic form in $(x,u)$,
condition \rf{CAL9} means that $\mu(x,u)=|Cx+Du|^2$
for some $C\in\CC^{m,n}$ and $D\in\CC^{m,m}$ such that
$D$ is not singular. In other words, representation
\rf{k2} takes place, and the inequality in \rf{k3} is
satisfied for $\lambda=0$.

To show that  the inequality in \rf{k3} is
satisfied for  $0<|\lambda|<1$, note that otherwise
there exist $p\in\CC^n$, $q\in\CC^m$,
 and $\xi\in\CC$ such that
\[  \AR{c}{p\\ q}'\AR{cc}{\xi A-I&\xi B\\ C & D}=0,\ \ 
\AR{c}{p\\ q}\neq0,\ \ |\xi|\in(0,1).\]
Then $q\neq0$ (otherwise $\xi p'A=p$, $p'B=0$, $p\neq0$
and hence the pair $(A,B)$ is not
stabilizable), and therefore it is possible to re-scale $(p,q)$
in such a way that $|q|=1$. Hence, for a 
solution $x,u$ of \rf{k7}
\begin{eqnarray*}
|Cx(t)+Du(t)|^2&\ge&|q'Cx(t)+q'Du(t)|^2\\
&=&|p'x(t)-\xi p'Ax(t)-\xi p'Bu(t)|^2\\
&=&|p'x(t)-\xi p'x(t+1)|^2,
\end{eqnarray*}
which implies that
\EQ{CAL10}{  
\sum_{t=0}^\infty|Cx(t)+Du(t)|^2\ge(1-|\xi|^2)|p'a|^2,}
contradicting the construdtion of $C,D$, which guarantees that
the maximal lower bound of the left side in \rf{CAL10} is zero
for all $a\in\CC^n$.
\end{itemize}

\subsubsection{Theorem~\ref{thm:kypstabdt}, 
the Case of Real Coefficients}
When the matrices $A,B,Q$ in (b) are real,
for every solution
$(x,u)$ of \rf{k7A} with $x(0)=a$ the conjugated pair
$(\bar x,\bar u)$ is a solution of \rf{k7A} with 
$x(0)=\bar a$ and the same value
of $\Phi$. Hence $V(\bar a)=V(a)$
for every $a\in\CC^n$, which implies that the (generally
complex) matrix $P=P'$ in the representation
$V(a)=-a'Pa$ actually has real coefficients.
Since, in this case, the Hermitian form
$\s(x,u)-V(x)+V(Ax+Bu)$ has real coefficients,
the matrices $C,D$ can also be chosen to be real.

\subsubsection{Proof of Theorem~\ref{thm:lqrdt}}
By \rf{k2} we have 
\[ \Phi=-a'Pa+\sum_{t=0}^\infty|Cx(t)+Du(t)|^2,\]
and it was already shown in the proof of 
Theorem~\ref{thm:kypstabdt},(b)$\Rightarrow$(a) that
$\inf\Phi=-a'Pa$. Hence $\Phi$ converges 
to its maximal lower bound if and only if the sum of
squares of $Cx+Du$ converges to zero.

\subsubsection{Proof of Theorem~\ref{thm:kyplmisdt}}
The implication (a)$\Rightarrow$(b) is trivial, as substituting
a non-zero pair $(x,u)$ from $\cl L(z)$ with $|z|=1$ into \rf{k5}
yields $\s(x,u)=\s_P(x,u)>0$.

To prove that (b) implies (a), assume that (b) is true but
(a) is not, which means that $0$ is not in the convex set
\[  \W=\{Q+E_0'PE_0-E_1'PE_1-S: \ \ S=S'>0,\ P=P'\},\]
where
\[ E_0=[I_n\ \ \ 0],\ \ \ E_1=[A\ \ \ B].\]
According to the Hahn-Banach theorem there exists a hyperplane
which separates (non-strictly) $\W$ from zero, i.e.
there exists matrix $H=H'\neq0$ such that
\EQ{CAL11}{\tr(XH)\le0\ \ \ \forall\ X\in\W.}
Using \rf{CAL11} with $X=Q-tI$ where $t\to0$ yields
$\tr(QH)\le0$.
Using \rf{CAL11} with $X=Q-I-tpp'$ where $t\to+\infty$
yields $\tr(Hpp')\ge0$ for every $p\in\CC^{n+m}$, i.e.
$H\ge0$. Similarly, using \rf{CAL11} with
$X=Q-I+t(E_0'PE_0-E_1'PE_1)$ where $t\to\pm\infty$
yields
$\tr(H(E_0'PE_0-E_1'PE_1))=0$ for every $P=P'$, i.e.
$E_0HE_0'=E_1HE_1'$. The last equality implies existence of
a unitary matrix $U$ such that $UH^{1/2}E_0'=H^{1/2}E_1'$,
or, equivalently, $E_0H^{1/2}U'=E_1H^{1/2}$. Let
$w_1,\dots,w_{n+m}$ be an orthonormal basis of
eigenvectors of $U'$, with $z_i\in\TT$ being the corresponding
eigenvalues. Define $x_i\in\CC^n$, $u_i\in\CC^m$ by
\[  e_i=\AR{c}{x_i\\ u_i}=H^{1/2}w_i.\]
By construction, $(x_i,u_i)\in\cl L(z_i)$, and hence
by assumption (b)
$e_i'Qe_i>0$ whenever $e_i\neq0$. On the other hand
\[ 0\ge\tr(QH)=
\tr(Q\sum_{i=1}^{n+m} e_ie_i')
=\sum_{i=1}^{n+m} e_i'Qe_i,\]
which means $e_i'Qe_i=0$ for all $i$. Hence $e_i=0$ for all $i$
and therefore $H=0$, which contradicts the construction.

To complete the proof, consider the case when $A,B,Q$ have
real coefficients. Then for every $P=P'$ such that
$\s_P>0$ we also have $\s_{\bar P}>0$, and hence,
for $\tilde P=0.5(P+\bar P)$,
\[  \s_{\tilde P}=0.5(\s_P+\s_{\bar P})>0.\]

\subsubsection{Proof of Theorem~\ref{thm:kyplminsdt}}
The implication (a)$\Rightarrow$(b) follows in the standard was
by substituting
an arbitrary pair $(x,u)$ from $\cl L(z)$ with $|z|=1$ into \rf{k5},
which
yields $\s(x,u)=\s_P(x,u)\ge0$.

To prove that (b) implies (a), consider the optimization task
\rf{k7}-\rf{k7B}, take any $K$ such that $A+BK$ is a Schur
matrix, and consider the Fourier transform repesentation
of $\Phi$ given by \rf{RIC7}. Since $\a(z)\ge0$ for all $z\in\TT$,
we have $\inf\Phi>-\infty$ for $a=0$. Therefore
$\inf\Phi>-\infty$ for every $a\in\CC^n$ which is reachable from
$x(0)=0$ in system \rf{k7}.
Since the pair $(A,B)$ is
controllable, we conclude that $\inf\Phi>-\infty$ for every $a\in\CC^n$.

We now use the same arguments as in the proof of 
Theorem~\ref{thm:kypstabdt} to establish that
$\inf\Phi=-a'Pa$ for some matrix $P=P'$ (real whenever
$A,B,Q$ are real). Finally, positive semidefiniteness of
$\s_P$ follows from the Bellman equation.

\subsubsection{KYP Proofs in Continuous Time}
In principle it is possible to translate all steps in the proofs of
Theorems \ref{thm:kypstabdt}-\ref{thm:kyplminsdt}
into a continuous-time format. However, there is a simple way
of deriving the CT versions from the DT ones.

Choose $r>0$ in such a way that the matrix $rI-A$ is not
singular. 
Let $\bar\CC=\CC\cup\{\infty\})$.
Consider the bijection $h_o:\ \bar\CC\mapsto\bar\CC$ 
and the linear bijection $h_1:\ \CC^n\times\CC^m\mapsto
\CC^n\times\CC^m$
which map $s\in\bar\CC$ to $z=h_0(s)\in\bar\CC$ 
and $(x,u)\in\CC^n\times\CC^m$ to
$(\tilde x,u)=h_1(x,u)\in\CC^n\times\CC^m$
according to
\[ \tilde x=\frac{rx-Ax-Bu}{\sqrt{2r}},\ \ \ 
z=\OR{ll}{\infty,& s=r,\\ -1,& s=\infty,\\
\frac{r+s}{r-s},& \tx{otherwise.}}\ \ 
\]
Define $\tilde A,\tilde B,\tilde\s$ by
\[ \tilde A=(rI+A)(rI-A)^{-1},\ \ 
\tilde B=\sqrt{2r}(rI-A)^{-1}B,\]
and
\[ \tilde\s(\tilde x,u)=\s(x,u)\ \ \tx{for}\ 
(\tilde x,u)=h_1(x,u)\]
(note that $\tilde A,\tilde B,\tilde\s$ will have real coefficients
whenever $A,B,\s$ have real coefficients).

Simple algebraic manipulations can be used to show that
\begin{itemize}
\item[(a)] for $(\tilde x,u)=h_1(x,u)$, equality
$sx=Ax+Bu$ is satisfied
if and only if $z\tilde x=\tilde A\tilde x+\tilde Bu$
(including the case $s=\infty$, $z=-1$, in which case
$sx=Ax+Bu$ is interpreted as $x=0$, as well as 
the case $z=\infty$, $s=r$, in which case
$z\tilde x=\tilde A\tilde x+\tilde Bu$ 
is interpreted as $\tilde x=0$);
\item[(b)] for every $(\tilde x,u)=h_1(x,u)$, the identity
\[ 2x'P(Ax+Bu)=(\tilde A\tilde x+\tilde Bu)'P(\tilde A\tilde x+\tilde Bu)
-\tilde x'P\tilde x\]
holds;
\item[(c)] $h_0(j\RR\cup\{\infty\})=\TT$;
\item[(d)] $s$ is an eigenvalue of $A$ if and only if
$z=h_0(s)$ is an eigenvalue of $\tilde A$.
\end{itemize}
In order to prove the CT statements for some $A,B,\s$, 
choose $r$ and construct $\tilde A,\tilde B,\tilde\s$ first.
For $z\in\TT\backslash\Lambda(\tilde A)$ 
define $\tilde\Pi=\tilde\Pi(z)$ by the identity
\[  u'\Pi(z)u=\tilde\s(\tilde x,u),\ \ \tx{subject to}\ \ 
z\tilde x=\tilde A\tilde x+\tilde Bu.\]
Then $\Pi(s)=\tilde\Pi(z)$ for $z=h_0(s)$, i.e. the
positive definiteness/semidefiniteness of $\tilde\Pi(z)$ on $\TT$
is determined by positive definitess/semidefiniteness of
$\Pi(s)$ for $s\in\RR\cup\{\infty\}$.
When matrices $C,D$ are given, define
$\tilde C,\tilde D$ by the identity
\[ Cx+Du=\tilde C\tilde x+\tilde Du\ \ \tx{for}\ \ 
(\tilde x,u)=h_1(x,u).\]
Now the DT statements of the KYP Lemma applied to 
$\tilde A,\tilde B,\tilde\s,\tilde\Pi$ (and, possibly, $\tilde C$,
$\tilde D$, $\tilde P=P$) prove the corresponding CT statements of the
KYP Lemma.

\subsection{Minimax Theorem Proofs}
This section contains the proofs of
the minimax theorems  associated with the KYP setup,
as well as the corresponding IQC statements.

\subsubsection{Proof of Theorems \ref{thm:minmaxdt}
and \ref{thm:minmaxct}}

The proof is based on associating the statements with the more
general setup of Theorem~\ref{thm:qminmax}.

In the DT case,  let $V=\ell^2_k$, $W=\ell^2_q$. The
functional $g$ defined by  \rf{k23},\rf{k25} is a quadratic form of
$(v,w,a)\in V\times W\times\RR^n$. Hence for
every fixed $a\in\RR^n$ it defines it defines unique quadratic
forms $\s,\mu$, bilinear form $p$, linear functions $f,h$, and
a constant $r$ such that representation \rf{CAL21} takes place.
According to Theorem \ref{thm:kypstabdt}, condition
$\Pi_{11}(z)\ge0$ (for $z\in\TT$), coupled with 
$\Pi_{11}(z_0)>0$ 
(both parts of  assumption \rf{k27}), implies that $g(v,0)$ has a finite lower bound.
Similarly, 
$\Pi_{22}(z)\le0$ (for $z\in\TT$), coupled with 
$\Pi_{11}(z_0)<0$, implies that $g(0,w)$ has a finite upper bound,
so condition (ii) of Theorem~\ref{thm:qminmax} is satisfied.
Finally, in terms of Fourier transforms we have
\[  \s(v)=\int_{\TT}\hat v'\Pi_{11}\hat vdm(z),\ 
\mu(w)=\int_{\TT}\hat w'\Pi_{22}\hat wdm(z),\ 
p(v,w)=\tx{Re}\int_{\TT}\hat v'\Pi_{12}\hat wdm(z).\]
Since  \rf{k27} implies that the matrix
\[  \int_{\TT}\AR{cc}{\hat v&0\\ 0&\hat w}'
\AR{cc}{\Pi_{11} & \e\Pi_{12}\\ \e\Pi_{21} & -\Pi_{22}}
\AR{cc}{\hat v&0\\ 0&\hat w}dm(z)=
\AR{cc}{ \int_{\TT}\hat v'\Pi_{11}\hat vdm(z) &
\e\int_{\TT}\hat v'\Pi_{12}\hat wdm(z)\\
\e\int_{\TT}\hat w'\Pi_{21}\hat vdm(z) &
-\int_{\TT}\hat w'\Pi_{22}\hat wdm(z)}\]
is positive semidefinite for all $v\in\ell^2_k$, $w\in\ell^2_q$,
condition (i) is satisfied with $c=\e^{-1}$.

According to  Theorem~\ref{thm:qminmax}
this means that the minimax equality holds. The bounds 
in (b) are now established as well. Since the partial optimal values
are quadratic functionals, the bounds establish their continuity.

The proofs for the CT case follow the same pattern, with CT
Fourier transform replacing the DT version.

\subsubsection{Proof of Theorems \ref{thm:minmaxiqcdt}
and \ref{thm:minmaxiqcct}}

Consider the DT case first. For every $(v_0,w_0)\in\D$, $x_0\in X_0$,
and $T>0$ consider functional $g$ from Theorem~\ref{thm:minmaxdt}
defined with $a=x(T+1)$. According to (b), 
for every fixed $w\in\ell^2_q$ the maximal lower bound
of $g(v,w)$ with respect $v\in\ell^2_k$ is not positive, i.e.
\[\sup_w\inf_v g(v,w)\le0.\] 

By
Theorem~\ref{thm:minmaxdt}, 
\[\inf_v\sup_w g(v,w)=\sup_w\inf_v g(v,w)\le0,\] 
which means that there exists a
sequence of signals $\{\tilde v_i\}_{i=1}^\infty\subset\ell^2_k$  
such that
$g(\tilde v_i,w)<1/i$ for all $i\in\{1,2,\dots\}$ and all $w\in\ell^2_q$.
In addition, Theorem~\ref{thm:minmaxdt} also claims that
$\sup_w g(v,w)$ is a continuous in the metric of the Hilbert space
defined by the associated quadratic form $\s$. 
Since, for $v\in\ell^2_k$,
\[  \s(v)=\int_{\TT}\hat v(z)'\Pi_{11}(z) v(z)dm(z),\]
and $\Pi_{11}$ is uniformly bounded on $\TT$, the norm
$\s(v)^{1/2}$ is majorated by the standard Hilbert space norm of
$\ell^2_k$, and hence 
$\sup_w g(v,w)$ is a continuous in the standard metric of $\ell^2_k$.
Accordingly, there exist a
sequence of signals $\{\tilde v_i\}_{i=1}^\infty\subset\ell^2_k$  
and a sequence of positive numbers $\d_i>0$
such that
$g(u,w)<1/i$ for all $i\in\{1,2,\dots\}$, $u\in\ell^2_k$,
and $w\in\ell^2_q$ such that $|u-\tilde v_i|<\d_i$.

For every $i\in\{1,2,\dots\}$ let
\[  v_*(t)=\OR{ll}{v(t),& t\le T,\\ \tilde v_i(t-T-1), & t>T.}\]
Since $\D$ is assumed to be weakly causally stable, 
there exist $(\tilde v,\tilde w)\in\D\cap(\ell^2_k\times\ell^2_q)$
such that
$\tilde v(t)=v(t)$,  $\tilde w(t)=w(t)$ for $t\le T$, and
$|\tilde v-v_*|<\d_i$. Due to the way in which $\tilde v_i,\d_i$
were chosen, for the corresponding solution
$\tilde x$ of 
\[ \tilde x(t+1)=A\tilde x(t)+B_1\tilde v(t)+B_2\tilde w(t),\ \ 
\tilde x(0)=x_0\]
we have
\[  \sum_{t>T}\s(\tilde x(t),\tilde v(t),\tilde w(t))<\frac1i.\]
Since  
\[  \sum_{t=0}^\infty\s(\tilde x(t),\tilde v(t),\tilde w(t))
\ge\kappa(x_0,v(0),w(0))\]
by the conditional IQC assumption, and
$x=\tilde x$, $v=\tilde v$, $w=\tilde w$ for $t\le T$ for all $i$,
we conclude (by letting $i\to\infty$) that
\[  \sum_{t=0}^T\s(x(t),v(t),w(t))
\ge\kappa(x_0,v(0),w(0)),\]
which proves the complete IQC.

The derivation in the CT time case follows the same steps,
with the definitions of $a$ and $v_*$ being modified to
$a=x(T)$ and
\[  v_*(t)=\OR{ll}{v(t),& t\le T,\\ \tilde v_i(t-T), & t>T.}\]

\section{Ackwnoledgements}
This paper was written as a technical addendum to the
Control Handbook article
"Integral Quadratic Constraints" by 
Alexandre Megretski, Ulf T. J\"onsson, Chung-Yao Kao,
and Anders Rantzer. In particular, it establishes the
claims of Theorem 2.1, and expands on the "minimax" approach
to post-feasibility analysis mentioned there.

The author is very grateful to Ulf J\"onsson and Anders Rantzer
for their encouragement and many helpful comments.

\end{document}